\documentclass[12pt]{article}
\usepackage{amsmath}
\usepackage{amssymb}
\usepackage{latexsym}
\usepackage{amsthm}
\usepackage{mathrsfs}
\usepackage{enumitem}
\usepackage[colorlinks,
            linkcolor=red,
            anchorcolor=blue,
            citecolor=green
            ]{hyperref}
\usepackage{graphicx}
\usepackage{diagbox}
\usepackage{tikz}
\usepackage{color}

\newtheorem{thm}{Theorem}[section]
\newtheorem{lem}[thm]{Lemma}

\newtheorem{cor}[thm]{Corollary}

\DeclareMathOperator{\red}{red}

\makeatother

\begin{document}

\begin{center}
{\large \bf  On pattern avoiding indecomposable permutations}
\end{center}

\begin{center}
Alice L.L. Gao$^{1}$,
Sergey Kitaev$^{2}$, and Philip B. Zhang$^{3}$\\[6pt]

$^{1}$Center for Combinatorics, LPMC-TJKLC\\
Nankai University, Tianjin 300071, P. R. China\\[6pt]

$^{2}$ Department of Computer and Information Sciences \\
University of Strathclyde, 26 Richmond Street, Glasgow G1 1XH, UK\\[6pt]

$^{3}$ College of Mathematical Science \\
Tianjin Normal University, Tianjin  300387, P. R. China\\[6pt]

Email: $^{1}${\tt gaolulublue@mail.nankai.edu.cn},
	   $^{2}${\tt sergey.kitaev@cis.strath.ac.uk},
           $^{3}${\tt zhangbiaonk@163.com}
\end{center}

\noindent\textbf{Abstract.}
Comtet introduced the notion of indecomposable permutations in 1972. A permutation is indecomposable if and only if it has no proper prefix which is itself a permutation. Indecomposable permutations were studied in the literature in various contexts. In particular, this notion has been proven  to be useful in obtaining non-trivial enumeration and equidistribution results on permutations.

In this paper, we give a complete classification of indecomposable permutations avoiding a classical pattern of length 3 or 4, and of indecomposable permutations avoiding a non-consecutive vincular pattern of length 3. Further, we provide a recursive formula for enumerating $12\cdots k$-avoiding indecomposable permutations for $k\geq 3$. Several of our results involve the descent statistic. We also provide a bijective proof  of a fact relevant to our studies.\\



\section{Introduction}\label{intro}
Let $[n]=\{1,\ldots,n\}$ and $\mathfrak{S}_n$ be the set of permutations of $[n]$.
Given $\pi=\pi_1\cdots \pi_n \in \mathfrak{S}_n$, let $i_{\pi}$ denote the smallest index such that $\pi_1\cdots \pi_{i_{\pi}}$ is a permutation of $[i_{\pi}]$. If $i_{\pi}=n$ then $\pi$  is \emph{indecomposable}; otherwise, $\pi$ is {\em decomposable}. For example, 23514 is indecomposable, while 31254 is decomposable.

For a permutation $\pi$ of a set $\{a_1,\ldots,a_n\}$, the {\em reduced form} of $\pi$, denoted  $\red(\pi)$, is the permutation of $[n]$ obtained from $\pi$ by replacing the $i$-th smallest element with $i$. For example, $\red(2537)=1324$. For any permutation $\pi$, we have $\pi=\pi^{(1)}\pi^{(2)}\cdots \pi^{(k)}$ for some $k\geq 1$, where $\pi^{(i)}$ is a permutation such that $\red(\pi^{(i)})$ is indecomposable for all $1\leq i\leq k$. We say that a $\pi^{(i)}$ is a {\em component} of $\pi$. For example, the permutation $312465$ has components 312, 4 and 65.

%For any $\pi\in\mathfrak{S}_n$, let $\pi=\pi^{(1)}\pi^{(2)}$, where  $\pi^{(1)}$ is a permutation of $\{1,\ldots,i_{\pi}\}$
%and $\pi^{(2)}$ is a possibly empty permutation of $\{i_{\pi}+1,i_{\pi}+2,\ldots,n\}$.

A ({\em permutation}) {\em pattern} is a permutation $\tau=\tau_1\cdots\tau_k$. We say that a permutation $\pi=\pi_1\cdots\pi_n$ {\em contains an occurrence} of $\tau$ if there are $1\leq i_1<\cdots< i_k\leq n$ such that $\pi_{i_1}\cdots \pi_{i_k}$ is order-isomorphic to $\tau$, that is, if the reduced form of $\pi_{i_1}\cdots \pi_{i_k}$ is $\tau$. If $\pi$ does not contain an occurrence of $\tau$, we say that $\pi$ {\em avoids}~$\tau$. This type of patterns is referred to as ``classical patterns''.
 For instance, the permutation 315267 contains several occurrences of the pattern 123, for example, the subsequences 356 and 157, while this permutation avoids the pattern 321. A comprehensive introduction to the theory of patterns in permutations can be found in~\cite{Kitaev2011Patterns}.

Another type of patterns of interest to us is {\em vincular patterns}, also known as {\em generalized patterns} \cite{Babson2000Generalized,Steingrimsson2010Generalized}, in occurrences of which some of the elements may be required to be adjacent in a permutation. We underline elements of a given pattern to indicate the elements that must be adjacent in any occurrence of the pattern. For example, the permutation $\pi=136254$ contains four occurrences of the pattern $1\underline{32}$, namely, the subsequences 162, 154, 354 and 254 (in each of these occurrences, the elements in $\pi$ corresponding to 2 and 3 in the pattern stay next to each other). On the other hand, $\pi$ contains only one occurrence of the pattern  $\underline{132}$, namely, 254. If all elements in an occurrence of a pattern are required to stay next to each other, which is indicated by underlying all elements in the pattern, such a pattern is called a {\em consecutive pattern}. A classical statistic {\em descent} is just an occurrence of the pattern $\underline{21}$. Vincular patterns play an important role in the theory of patterns in permutations and words (see Sections 3.3 and 3.4 in \cite{Kitaev2011Patterns} for details).

The notion of indecomposable permutations (also known as {\em irreducible permutations} or {\em connected permutations}) was introduced by Comtet \cite{Comtet1972Sur,Comtet1974Advanced}. Comtet was the first one to show the not so difficult to see fact that the ordinary generating function for the number $I_n$ of indecomposable permutations of length $n$ is
$$
\sum_{n=1}^{\infty}I_n x^n =  1-\frac{1}{\sum_{k\ge 0}k!x^k}.
$$
These numbers  begin with 1, 1, 3, 13, 71, 461, 3447, 29093, $\ldots$ for $n\geq 1$, and appear as the sequence A003319 in the {\em On-Line Encyclopedia of Integer Sequences} ({\em OEIS}) \cite{oeis}.

Indecomposable permutations appear in various contexts in the literature, for example, see  \cite{Claesson2008Classification,CKS,Cori2009Hypermaps, Cori2009Indecomposable,King2006Generating}. In particular, in \cite[Section 4]{Claesson2008Classification}, indecomposable permutations are used to define a bijection between 231- and 321-avoiding permutations (finding various bijections essentially between these sets was the subject of several papers in the literature). Also, indecomposable pattern avoiding permutations are a key object in \cite{CKS} to find a bijection between permutations in question and so-called {\em $\beta(1,0)$-trees} that has as corollaries a number of equidistribution results on permutations, $\beta(1,0)$-trees and certain types of planar maps. Finally, we note that indecomposable pattern avoiding permutations were first studied by B\'ona in \cite{Bona1997Exact}, where essentially 2431-avoiding indecomposable permutations are enumerated (indecomposable permutations are defined up to reverse in \cite{Bona1997Exact}) and linked in a bijective way to labeled plane trees of a certain type. 

In this paper, we study interrelations (taking into account the descent statistic, which comes ``for free'') between pattern avoiding permutations and their indecomposable counterparts for classical patterns of length 3 and 4, vincular non-consecutive patterns of length 3, and the increasing classical pattern of arbitrary length (patterns of length 2 are trivial).  We use our results and known enumeration formulas for pattern avoiding permutations to enumerate indecomposable pattern avoiding permutations. Some of the obtained numbers appear in the OEIS suggesting a number of bijective questions. 

The paper is organized as follows. In Section~\ref{prelim}, we introduce generating functions to be studied in this paper and state known pattern avoiding results to be used. In Section~\ref{sec-class-pattern-avoid}, we study indecomposable permutations avoiding classical patterns of length 3 and 4, as well as the classical pattern $12\cdots k$ for $k\geq 3$.  In Section~\ref{vin-sec}, we study indecomposable permutations avoiding a vincular non-consecutive pattern of length 3, in particular, presenting a bijective result in Theorem~\ref{thm-bij}.  Finally, in Section~\ref{research}, we discuss directions of further research.

\section{Preliminaries}\label{prelim}

Let $A_{n}^{\sigma}$ and $I_{n}^{\sigma}$ be the number of $\sigma$-avoiding permutations of $[n]$ and $\sigma$-avoiding indecomposable permutations of $[n]$, respectively.
For $0\leq i \leq n-1$, let $A_{n,i}^{\sigma}$ and $I_{n,i}^{\sigma}$ be the number of $\sigma$-avoiding permutations of $[n]$ and $\sigma$-avoiding indecomposable permutations of $[n]$ with $i$ descents, respectively.
Thus, for $n\geq 1$, we have
$$A_{n}^{\sigma}=\sum_{i=0}^{n-1}A_{n,i}^{\sigma}~~~~\mathrm{and}~~~~
I_{n}^{\sigma}=\sum_{i=0}^{n-1}I_{n,i}^{\sigma}.$$
Let $A^{\sigma}(x)$, $A^{\sigma}(x,q)$, $I^{\sigma}(x)$ and $I^{\sigma}(x,q)$ be the generating functions for $A_{n}^{\sigma}$, $A_{n,i}^{\sigma}$, $I_{n}^{\sigma}$ and $I_{n,i}^{\sigma}$, respectively. That is,
$$A^{\sigma}(x,q)=\sum_{n=0}^{\infty}\sum_{i=0}^{n-1}A_{n,i}^{\sigma}x^nq^i,~~~I^{\sigma}(x,q)=\sum_{n=1}^{\infty}\sum_{i=0}^{n-1}I_{n,i}^{\sigma}x^nq^i,$$
$$A^{\sigma}(x)=A^{\sigma}(x,1)=\sum_{n=0}^{\infty}A_{n}^{\sigma}x^n~~\mbox{and}~~I^{\sigma}(x)=I^{\sigma}(x,1)=\sum_{n=1}^{\infty}I_{n}^{\sigma}x^n.$$
If for patterns $\sigma_1$ and $\sigma_2$, $A_{n}^{\sigma_1}=A_{n}^{\sigma_1}$ for all $n\geq 0$ then $\sigma_1$ and $\sigma_2$ are {\em Wilf-equivalent}.

For a permutation $\pi=\pi_1\cdots\pi_n$, its {\em reverse} is the permutation $r(\pi)=\pi_n\pi_{n-1}\cdots \pi_1$ and its {\em complement} is the permutation $c(\pi)=(n+1-\pi_1)(n+1-\pi_2)\cdots(n+1-\pi_n)$. For example, if $\pi=32145$ then $r(\pi)=54123$ and $c(\pi)=34521$.  The reverse and complement operations are called {\em trivial bijections}. It is easy to see that for any pattern $\sigma$, this pattern is Wilf-equivalent to $r(\sigma)$ and $c(\sigma)$. Another useful property of trivial bijections is that their composition preserves the property of being decomposable (and thus the property of being indecomposable), which is easy to see.

One of concerns in this paper is to find interrelations between  $I^{\sigma}(x,q)$ and $A^{\sigma}(x,q)$ for certain $\sigma$s. We note that throughout this paper we implicitly use the fact that an occurrence of a descent cannot start in one component of a permutation and end in another one.

In the rest of this section we review a number of permutation pattern avoidance results relevant to this paper. These results come from \cite[Chapters 6 and 7]{Kitaev2011Patterns}, where references to the original sources can be found.

\begin{lem}\label{catalan}
For $\sigma \in \mathfrak{S}_3$ and $n\geq 0$, $A^{\sigma}_{n}=C_n=\frac{1}{n+1}\binom{2n}{n}$, the $n$-th Catalan number.
Thus, $$A^{\sigma}(x)=C(x)=\frac{1-\sqrt{1-4 x}}{2x},$$
the generating function for the Catalan numbers satisfying  $xC(x)^2-C(x)+1=0$. For $n\geq 0$, the Catalan numbers  begin with $1,1,2,5,14,42,\ldots$, which is the sequence A000108 in the OEIS  \cite{oeis}.
\end{lem}

Next result links the well-known {\em Bell numbers} to pattern avoiding permutations. The Bell numbers begin with $1, 1, 2, 5, 15, 52, 203,\ldots$ for $n\geq 0$, and this is the sequence A000110 in the OEIS  \cite{oeis}.

\begin{lem}\label{bell}
When $\sigma\in\{1\underline{23},1\underline{32},3\underline{12},3\underline{21}\}$,  for $n\geq 0$, we have $$A^{\sigma}_n=B_n,$$
where $B_n$ is the $n$-th Bell number, which is the number of set partitions of $[n]$.

When
$\sigma\in\{2\underline{13},2\underline{31}$\}, for $n\geq 0$, we have
$$A^{\sigma}_n=C_n.$$
\end{lem}

We next turn our attention to classical patterns of length 4. Table~\ref{tab} presents three Wilf-equivalence classes in this case.

\begin{table}[!htb]
 \renewcommand{\arraystretch}{1.3}
\begin{center}
\begin{tabular}{|c|l|}
\hline
1 & 1234,4321,~~1243,2134,3421,4312,~~1432,2341,3214,4123,~~2143,3412 \\
\hline
2&1342,2431,3124,4213,~~1423,2314,3241,4132,~~2413,3142\\
\hline
3&1324,4231\\
\hline
\end{tabular}
\end{center}
\caption{The three Wilf-equivalence classes for pattern avoidance of length 4. Spaces on a line are used to indicate patterns equivalent via trivial bijections.}\label{tab}
\end{table}

$A^{1234}_n$ for $n\geq 0$ begins with $1, 1, 2, 6, 23, 103,\ldots$ (this is the sequence \cite[A005802]{oeis}), and we have the following lemma.

\begin{lem}\label{ge-1234}
For $\sigma\in \{1234,4321,1243,2134,3421,4312,1432,2341,3214,4123,2143,3412\}$, we have
$$E(x):=A^{\sigma}(x)= \frac{1}{6 x^2} \left( 1+5 x-(1-9 x)^{\frac{3}{4}} (1-x)^{\frac{1}{4}} \, _2F_1\left(-\frac{1}{4},\frac{3}{4};1;\frac{64 x}{(x-1) (1-9 x)^3}\right) \right).$$
Moreover, for $n\geq 1$, we have
$$E_n:=A^{\sigma}_n= 2  \sum_{k=0}^{n} \binom{2 k}{k} \binom{n}{k}^2  \frac{3k^2+2k-2kn-n+1}{(k+1)^2 (k+2) (n-k+1)}.$$
\end{lem}

An exact enumeration for 1342-avoiding permutations and the corresponding generating function are given by B{\'o}na \cite{Bona1997Exact}.  The corresponding sequence for $n\geq 0$ begins with $1, 1, 2, 6, 23, 103, 512, 2740,\ldots$ (A022558 in \cite{oeis}) and the following lemma holds.

\begin{lem}\label{ge-1342}
For $\sigma\in \{1342,2431,3124,4213,1423,2314,3241,4132,2413,3142\}$,
we have
$$F(x):=A^{\sigma}(x)=\frac{32 x}{1+20 x-8 x^2-(1-8 x)^{3/2}}=\frac{1+20 x-8 x^2}{2 (1+x)^3}+\frac{(1-8 x)^{3/2}}{2 (1+x)^3}.$$
Moreover, for $n\geq 1$, we have
$$F_n:=A^{\sigma}_n=\frac{7n^2-3n-2}{2} (-1)^{n-1} + 3 \sum_{i=2}^{n} (-1)^{n-i} 2^{i+1} \frac{(2i-4)!}{i!(i-2)!}  \binom{n-i+2}{2}.$$
\end{lem}

However, no formula for $A^{1324}_{n}$ is known, only a recurrence relation is discovered \cite{Marinov2002/03Counting},  and an algorithm for counting the number of 1324-avoiding permutations was given in \cite{Conway20151324,Johansson2014Using}.
 For recent developments on the bounds, see  \cite{Bevan2015Permutations,Bona2015new}.
The corresponding sequence for $n\geq 0$ begins with $1, 1, 2, 6, 23, 103, 513, 2762,\ldots$; see A061552 in \cite{oeis}.

\section{Indecomposable permutations avoiding classical patterns}\label{sec-class-pattern-avoid}

Patterns are permutations, and we distinguish two cases according to whether or not they are decomposable. We start with an easier case.

\subsection{Indecomposable patterns}

Here we deal with the following patterns:

231, 312, 321,

2341, 2413, 2431, 3142, 3241, 3412, 3421, 4123, 4321, 4132, 4213, 4231, 4312.
%
%The generating functions $A^{\sigma}(x)$ and $I^{\sigma}(x)$ have the following close connection:
%$$A^{\sigma}(x)= \frac{1}{1-I^{\sigma}(x)},$$
%equivalently,

We first establish a property holding for any indecomposable pattern $\sigma$.

\begin{lem}\label{inde}
If $\sigma$ is an indecomposable pattern, then $I^{\sigma}(x,q)$ satisfies
$$I^{\sigma}(x,q)= 1-\frac{1}{A^{\sigma}(x,q)}.$$
\end{lem}

\proof
For any permutation $\pi$, an occurrence of $\sigma$ cannot start in one component and end in another one, which would contradict $\sigma$ being irreducible. Similarly, a descent cannot start in one component and end in another one. Hence, the generating function for $\sigma$-avoiding permutations with $k$ components is $[I^{\sigma}(x,q)]^k$ and
$$A(x,q)=1+\sum_{k=1}^{\infty}\left( I^{\sigma}(x,q) \right)^k=\frac{1}{1-I^{\sigma}(x,q)},$$
where ``1+'' corresponds to the empty permutation. This gives the desired result.
\qed

Combining Lemma~\ref{inde} ($q=1$) and Lemma~\ref{catalan} we obtain the following theorem, which can also be derived, e.g. from considerations in \cite{Claesson2008Classification}.

\begin{thm}\label{inde-length-3}
For $\sigma\in\{312,321,231\}$, we have
\begin{align}\label{eq-catalan}
I^{\sigma}(x)=\frac{1-\sqrt{1-4x}}{2}.
\end{align}
Thus, for $n\geq 1$, $I^{\sigma}_n=C_{n-1}$, the $(n-1)$-th Catalan number.
\end{thm}

Combining Lemma~\ref{inde} ($q=1$) and Lemma~\ref{ge-1342} we obtain the following theorem essentially established in~\cite{Bona1997Exact}.

\begin{thm}
For $\sigma\in \{2431,4213,3241,4132,2413,3142\}$,
we have
$$I^{\sigma}(x)=\frac{-1+12 x+8 x^2+(1-8 x)^{3/2}}{32x}.$$ The initial values for $I^{\sigma}_n$ in this case are $1, 1, 3, 12, 56, 288, 1584, 9152,\ldots$ for $n\geq 1$, and this is the sequence A000257 in the OEIS~\cite{oeis}.
\end{thm}

Similarly, one can combine Lemma~\ref{inde} ($q=1$) and Lemma~\ref{ge-1234} to obtain a formula for $I^{\sigma}(x)$, where $\sigma\in \{4321,3421,4312,2341,4123,3412\}$. The initial values for $I^{\sigma}_n$ in this case are $1, 1, 3, 12, 56, 289, 1603, 9391,\ldots$ for $n\geq 1$, and this sequence is not in the OEIS~\cite{oeis}.

However, we cannot obtain a formula for $I^{4231}(x)$ using Lemma~\ref{inde} because no formula is known for $A^{4231}(x)$. The initial values for $I^{4231}_n$ are $1, 1, 3, 12, 56, 289, 1604, 9415,\ldots$ for $n\geq 1$, and this sequence is not in the OEIS~\cite{oeis}.

\subsection{Decomposable patterns}

The patterns here we deal with are

123, 132, 213,

1234, 1243, 1324, 1342, 1423, 1432, 2134, 2143, 2314, 3124, 3214.

\subsubsection{Decomposable patterns of length 3}

\noindent
{\bf Pattern 123.} We first give a description of $123$-avoiding decomposable permutations.
\begin{lem}\label{123-decom}
Let $\pi=\pi^{(1)}\pi^{\prime}\in \mathfrak{S}_n$, where $\pi^{(1)}$ is the component in $\pi$ formed by the elements in $\{1,\ldots,i_{\pi}\}$. Then $\pi$ is a $123$-avoiding decomposable permutation if and only if
 $$\pi^{(1)}= i_{\pi} (i_{\pi}-1) \cdots 1, \pi^{\prime}=n(n-1)\cdots  (i_{\pi}+1)\mbox{ and }1\leq i_{\pi}\leq n-1.$$
\end{lem}

\proof
The backward direction is straightforward to see since no occurrence of the pattern $123$ can start in $\pi^{(1)}$.

For the forward direction, since $\pi=\pi_1\cdots\pi_n$ is decomposable, one must have $1\leq i_{\pi}\leq n-1$.  If $\pi_i<\pi_j$ for some $1\leq i<j\leq i_{\pi}$ then $\pi_i\pi_j\pi_n$ is an occurrence of the pattern 123; contradiction. If $\pi_i<\pi_j$ for some $i_{\pi}+1\leq i<j\leq n$ then $\pi_1\pi_i\pi_j$ is an occurrence of the pattern 123; contradiction. Thus, we obtain the desired result. \qed

Next we derive a relation between $I^{123}(x,q)$ and $A^{123}(x,q)$, which will give formulas for $I^{123}(x)$  and $I^{123}_{n}$. The initial values for $I^{123}_{n}$ for $n\ge 1$ begin with $1, 1, 3, 11, 38, 127, 423, \ldots$. This sequence does not appear in the OEIS \cite{oeis}.

\begin{thm}\label{thm-123}
We have that
\begin{align}\label{123-1}
A^{123}(x,q)=I^{123}(x,q)+\frac{x^2}{(1-xq)^2}+1,
\end{align}
 \begin{align}\label{123-11}
 I^{123}(x) = \frac{1-\sqrt{1-4 x}}{2x}- \frac{x^2}{(1-x)^2}-1,
 \end{align}
and for $n\geq 1$,
 \begin{align}\label{123-111}
I^{123}_{n} = C_{n}- (n-1).
 \end{align}
\end{thm}

\proof
By Lemma \ref{123-decom}, the generating function for all 123-avoiding decomposable permutations is
$\sum_{i\geq 1}x^iq^{i-1}\sum_{j\geq 1}x^jq^{j-1}=\frac{x^2}{(1-xq)^2}$. Since any 123-avoiding permutation is either indecomposable, or decomposable, or the empty permutation, we obtain the relation given by \eqref{123-1}. Letting $q=1$ in \eqref{123-1} and using Lemma~\ref{catalan}, we obtain \eqref{123-11}. Finally, by Lemma~\ref{catalan} and the fact that $\frac{x^2}{(1-x)^2}=\sum_{n\geq 1}(n-1)x^n$, we obtain \eqref{123-111}.
This completes the proof.
\qed

We note that (\ref{123-111}) appears in Proposition 9 in \cite{Disanto}. \\

\noindent
{\bf Patterns 132 and 213.} We begin with a description of $132$-avoiding decomposable permutations.
\begin{lem}\label{132-decom}
Let $\pi=\pi^{(1)}\pi^{\prime}\in \mathfrak{S}_n$, where $\pi^{(1)}$ is the component in $\pi$ formed by the elements in $\{1,\ldots,i_{\pi}\}$. Then $\pi$ is a $132$-avoiding decomposable permutation if and only if $\pi^{(1)}$ is a $132$-avoiding decomposable permutation,  $1\leq i_{\pi}\leq n-1$, and $\pi^{\prime}=(i_{\pi}+1)(i_{\pi}+2)\cdots n$.
\end{lem}

\proof
The backward direction is easy to see since an occurrence of the pattern $132$ cannot start in $\pi^{(1)}$ in this case.

For the forward direction, since $\pi$ is decomposable, we have $1\leq i_{\pi}\leq n-1$. Moreover, since $\pi$ is $132$-avoiding, $\pi^{(1)}$ must be $132$-avoiding. Finally, if $\pi_i>\pi_j$ for $i_{\pi}+1\leq i<j\leq n$ then $\pi_1\pi_i\pi_j$ is an occurrence of the pattern 132; contradiction.
\qed

Next we find a relation between $A^{\sigma}(x,q)$ and $I^{\sigma}(x,q)$ for $\sigma\in\{132,213\}$, which will give us formulas for $I^{\sigma}(x)$ and $I^{\sigma}_n$. The initial values for $I^{132}_n=I^{213}_n$ for $n\ge 1$ begin with $1, 1, 3, 9, 28, 90, 297, 1001,\ldots$, which is the sequence A000245 in \cite{oeis}.

\begin{thm}\label{thm-indecom-132}
For $\sigma\in\{132,213\}$, we have that
\begin{align}\label{eq-132}
A^{\sigma}(x,q)&=1+\frac{1}{1-x}I^{\sigma}(x,q),
\end{align}
\begin{align}\label{132-11}
I^{\sigma}(x)&=(1-x) \left(\frac{1 - \sqrt{1-4 x} }{2 x}-1\right),
\end{align}
$I^{\sigma}_{1}=1$, and for $n\ge 2$,
\begin{align}\label{132-111}
I^{\sigma}_{n} = C_{n}- C_{n-1}.
\end{align}
\end{thm}
\proof

Let $\sigma=132$. By Lemma \ref{132-decom}, the generating function for all 132-avoiding decomposable permutations is
$I^{132}(x,q)\sum_{j\geq1}x^j=\frac{x}{1-x}I^{132}(x,q)$.
Similarly to the proof of Theorem \ref{thm-123}, we have
\begin{align*}
A^{132}(x,q)&=1+I^{132}(x)+\frac{x}{1-x}I^{132}(x,q)\\
&=1+\frac{1}{1-x}I^{132}(x,q)
\end{align*}
giving \eqref{eq-132} for $\sigma=132$. Letting $q=1$ in \eqref{eq-132} and using Lemma \ref{catalan}, we obtain
\begin{align*}
I^{132}(x)&=(1-x)(A^{132}(x)-1)\\
&=(1-x) \left(\frac{1 - \sqrt{1-4 x} }{2 x}-1\right)
\end{align*}
giving \eqref{132-11} for $\sigma=132$. From \eqref{132-11}, $I^{132}_{1}=1$ and, for $n\ge 2$ \eqref{132-111} follows for $\sigma=132$.

Finally, since the composition of reverse and complement preserves the property of being irreducible, and this composition applied to 132 gives 213, we have that \eqref{eq-132}, \eqref{132-11} and \eqref{132-111} hold for $\sigma = 213$. \qed

Note that \eqref{132-111} follows directly from Lemma~\ref{132-decom}. Indeed, in a decomposable  $132$-avoiding permutation, the largest element $n$ must be the rightmost element, and the number of such permutations is $C_{n-1}$, while the number of all 132-avoiding permutations of length $n$ is $C_n$. Also, note that (\ref{132-111}) appears in Proposition 9 in \cite{Disanto}.

\subsubsection{Decomposable patterns of length 4}

Recall that applying the composition of reverse and complement, indecomposable permutations stay indecomposable, while applying that composition to the patterns we see that  each of  2143, 1324 and 1234  goes to itself, while for the remaining eight decomposable patterns this operation gives

$$I^{2314}(x,q)=I^{1423}(x,q),\ I^{3124}(x,q)=I^{1342}(x,q),$$
$$I^{3214}(x,q)=I^{1432}(x,q),\ \mbox{and } I^{2134}(x,q)=I^{1243}(x,q).$$
Thus, we only need to consider seven  decomposable patterns of length 4. \\

\noindent
{\bf Patterns 2314 and 3124.} We start with a description of $2314$-avoiding and 3124-avoiding decomposable permutations.

\begin{lem}\label{2314-decom}
Let $\pi=\pi^{(1)}\pi^{\prime}\in \mathfrak{S}_n$, where $\pi^{(1)}$ is a permutation of $\{1,\ldots,i_{\pi}\}$. Then $\pi$ is a $2314$-avoiding (resp., $3124$-avoiding) decomposable permutation if and only if
 $\pi^{(1)}$ is $231$-avoiding (resp., $312$-avoiding) and $\pi^{\prime}$ is $2314$-avoiding (resp., $3124$-avoiding).
\end{lem}

\proof
For the backward direction, because $\pi^{(1)}$ is $231$-avoiding (resp., $312$-avoiding) at most two elements in a possible occurrence of the pattern $2314$ (resp., $3124$) can be in $\pi^{(1)}$. But then $\pi^{\prime}$  contains an element smaller than an element in  $\pi^{(1)}$, which is impossible, and thus $\pi$ is $2314$-avoiding (resp., $3124$-avoiding).

For the forward direction, since $\pi=\pi_1\cdots\pi_n$ is decomposable, we have $1\leq i_{\pi}\leq n-1$. Also, clearly $\pi^{\prime}$ is $2314$-avoiding (resp., $3124$-avoiding). Now, if $\pi^{(1)}$ would contain an occurrence of the pattern $231$ (resp., $312$) then together with $\pi_n$ it would form an occurrence of the pattern $2314$ (resp., $3124$); contradiction. Thus $\pi^{(1)}$ is $231$-avoiding (resp., $312$-avoiding). \qed

Both of the sequences $I^{2314}_n$ and $I^{3124}_n$ for $n\geq 1$ begin with $1, 1, 3, 13, 65, 350, 1979,$ $11612,\ldots$, which does not appear in the OEIS \cite{oeis}.

\begin{thm}\label{thm-2314} We have
\begin{align}\label{2314-equ}
 A^{2314}(x,q) = 1+ I^{2314}(x,q) + I^{231}(x,q) \left(A^{2314}(x,q)-1\right)\mbox{ and}
\end{align}
\begin{align}\label{3124-equ}
 A^{3124}(x,q) = 1+ I^{3124}(x,q) + I^{312}(x,q) \left(A^{3124}(x,q)-1\right).
\end{align}
Further, for $\sigma\in\{2314,3124\}$, we have that
$$I^{\sigma}(x)=\frac{1}{2} \left(\sqrt{1-4 x}+1\right) \left(\frac{32 x}{1+20 x-8 x^2-(1-8 x)^{3/2}}-1\right),$$
$I^{\sigma}_{1}=1$, and for $n\geq 2$,
 $$
 I^{\sigma}_{n} =F_n- \sum_{i=0}^{n-2} C_{i}F_{n-1-i},
 $$
 where $F_n$ is defined in Lemma~\ref{ge-1342}.
\end{thm}

\proof
By Lemma \ref{2314-decom}, the generating function for $2314$-avoiding decomposable permutations is $I^{231}(x,q)(A^{2314}\left(x,q)-1\right)$, where ``$-1$'' corresponds to excluding the empty permutation as a possibility for $\pi^{\prime}$. Similarly, the generating function for $3124$-avoiding decomposable permutations is $I^{312}(x,q)(A^{3124}\left(x,q)-1\right)$.

Note that each $\sigma$-avoiding permutation is either the empty permutation, or an indecomposable permutation or a decomposable ones. This observation shows \eqref{2314-equ} and \eqref{3124-equ}.

Let $q=1$ in  \eqref{2314-equ}. Combing with Lemma \ref{ge-1342} and \eqref{eq-catalan}, we obtain
\begin{align}
 I^{2314}(x) & =  \left(A^{2314}(x)-1\right) \cdot \left(1- I^{231}(x)\right) \\
 & = \left(1-xC(x)\right) \left(F(x)-1\right)\label{2314-2}\\
 &=\frac{1}{2} \left(\sqrt{1-4 x}+1\right) \left(\frac{32 x}{1+20 x-8 x^2-(1-8 x)^{3/2}}-1\right).
\end{align}
 Hence, by \eqref{2314-2}, we have that $I^{2314}_1=1$ and for $n\geq 2$,
\begin{align*}
 I^{2314}_{n} &=F_n+C_{n-1}- \sum_{i=0}^{n-1} C_{i}F_{n-1-i}\\
 &=F_n- \sum_{i=0}^{n-2} C_{i}F_{n-1-i}.\label{2314-eqqq}
\end{align*}
It is straightforward to provide essentially the same derivations
 for the case of 3124-avoiding indecomposable permutations, which completes the proof. 
\qed

\ \\

\noindent
{\bf Pattern 3214.} We begin with a description  of $3214$-avoiding decomposable permutations. Our proof of next lemma is similar to the proof of Lemma~\ref{2314-decom} and thus is omitted.

\begin{lem}\label{3214-decom}
Let $\pi=\pi^{(1)}\pi^{\prime}\in \mathfrak{S}_n$, where $\pi^{(1)}$ is a permutation of $\{1,\ldots,i_{\pi}\}$. Then $\pi$ is a $3214$-avoiding decomposable permutation if and only if
$\pi^{(1)}$ is $321$-avoiding and $\pi^{\prime}$ is a $3214$-avoiding.
\end{lem}

The initial values $I^{3214}_n$ for $n\geq 1$ begin $1, 1, 3, 13, 65, 351, 1999, 11872, \ldots$ and this sequence is not in the OEIS \cite{oeis}.

\begin{thm} We have
\begin{equation}\label{3214-rel}
A^{3214}(x,q) = 1+ I^{3214}(x,q) + I^{321}(x,q) \left(A^{3214}(x,q)-1\right).
\end{equation}
Moreover, for $E(x)$ and $E_n$ defined in Lemma~\ref{ge-1234}, we have
$$I^{3214}(x)=\frac{1}{2} \left(\sqrt{1-4 x}+1\right) \left(E(x)-1\right),$$
$I^{3214}_1=1$ and for $n\geq 2$
 $$
 I^{3214}_{n} =E_n- \sum_{i=0}^{n-2} C_{i} E_{n-1-i}.
 $$
\end{thm}

\proof
We can proceed similarly to the proof of Theorem \ref{thm-2314} to prove \eqref{3214-rel}. Further, assuming that $q=1$ in \eqref{3214-rel}, one can apply Lemma \ref{ge-1234} and  \eqref{eq-catalan}, to obtain
 \begin{align*}
 I^{3214}(x) & =  \left(A^{3214}(x)-1\right) \cdot \left(1- I^{321}(x)\right) \\
 & = \left(1-xC(x)\right)\left(E(x)-1\right)\\
 &=\frac{1}{2} \left(\sqrt{1-4 x}+1\right) \left(E(x)-1\right).
 \end{align*}
From the last derivation, the formula for  $I^{3214}_{n}$ holds. %\ref{ge-1234},
\qed

\ \\

\noindent
{\bf Pattern 2143.} We begin with a description of $2143$-avoiding decomposable permutations.

\begin{lem}\label{2143-decom}
Let $\pi=\pi^{(1)}\pi^{\prime}\in \mathfrak{S}_n$, where $\pi^{(1)}$ is a permutation of $\{1,\ldots,i_{\pi}\}$. Then $\pi$ is a $2143$-avoiding decomposable permutation if and only if one of the following two conditions holds:
  \begin{itemize}
   \item $\pi^{(1)}=1$ and $\pi^{\prime}$ is $2143$-avoiding.
   \item $2\leq i_{\pi}\leq n-1$, $\pi^{(1)}$ is $2143$-avoiding and $\pi^{\prime}=(i_{\pi}+1)(i_{\pi}+2)\cdots n$ .
  \end{itemize}
\end{lem}

\proof
For the forward direction, since $\pi$ is decomposable, we have $1\leq i_{\pi}\leq n-1$.  There are two cases to consider:
\begin{itemize}
\item $i_{\pi}=1$. It is clear that $\pi^{(1)}=1$ in this case does not affect $\pi^{\prime}$, so  $\pi^{\prime}$ must be a $2143$-avoiding permutation of $\{2,\ldots,n\}$.
\item $2\leq i_{\pi}\leq n-1$. Since $\pi^{(1)}$ is indecomposable of length at least 2, there exist
$1\leq j_1<j_2\leq i_{\pi}$ such that $\pi_{j_1}>\pi_{j_2}$. But then to avoid an occurrence of the pattern $2143$ involving $\pi_{j_1}$ and $\pi_{j_2}$,  $\pi^{\prime}$ must be the increasing permutation $(i_{\pi}+1)(i_{\pi}+2)\cdots n$.
\end{itemize}
The backward direction is easy to see using similar considerations as above.\qed

The initial values $I^{2143}_n$ for $n\geq 1$ are  $1, 1, 3, 13, 63, 330, 1838, 10758,\ldots$, and this sequence is not in the OEIS \cite{oeis}.

\begin{thm} We have
\begin{equation}\label{eq-2142-case}
A^{2143}(x,q) = 1+ I^{2143}(x,q) + x \left(A^{2143}(x,q)-1\right) + \frac{x}{1-x} \left(I^{2143}(x,q)-x\right).
\end{equation}
Moreover, for $E(x)$ and $E_n$ defined in Lemma~\ref{ge-1234}, we have
$$ I^{2143}(x)  =  (1-x)^2 E(x) +2x-1,
$$
$I^{2143}_n=1$ and for $n\geq 2$,
$$I^{2143}_n=E_n-2E_{n-1}+E_{n-2}.$$
\end{thm}

\proof
By Lemma \ref{2143-decom}, \eqref{eq-2142-case} follows. Further, letting $q=1$ in \eqref{eq-2142-case} and using Lemma~\ref{ge-1234}, it follows that
\begin{align*}
 I^{2143}(x)&=(1-x)^2 A^{2143}(x) +2x-1\\
 &=(1-x)^2 E(x) +2x-1,
\end{align*}
from which the formula for $I^{2143}_n$ follows.  \qed

\ \\

\noindent
{\bf Pattern 2134.} We first give a description of $2134$-avoiding decomposable permutations.

\begin{lem}\label{2134-decom}
Let $\pi=\pi^{(1)}\pi^{\prime}\in \mathfrak{S}_n$, where $\pi^{(1)}$ is a permutation of $\{1,\ldots,i_{\pi}\}$. Then $\pi$ is a $2134$-avoiding decomposable permutation if and only if one of the following two conditions holds:
  \begin{itemize}
   \item $\pi^{(1)}=1$ and $\pi^{\prime}$ is $2134$-avoiding.
   \item $2\leq i_{\pi}\leq n-1$, $\pi^{(1)}$ is $213$-avoiding, and $\pi^{\prime}=n(n-1) \cdots (i_{\pi}+1)$.
  \end{itemize}
\end{lem}

\proof
For the forward direction, since $\pi=\pi_1\cdots\pi_n$ is decomposable, we have $1\leq i_{\pi}\leq n-1$.  There are two cases to consider:
\begin{itemize}
\item $i_{\pi}=1$. It is clear that $\pi^{(1)}$ does not affect the rest of the permutation, and $\pi^{\prime}$ must be a $2134$-avoiding permutation of $\{2,\ldots,n\}$.
\item $2\leq i_{\pi}\leq n-1$. Then $\pi^{(1)}$ must be a $213$-avoiding permutation, or else it would form an occurrence of the pattern  2134 with $\pi_n$.  Moreover, since $\pi^{(1)}$ is indecomposable, there exist
$1\leq j_1<j_2\leq i_{\pi}$ such that $\pi_{j_1}>\pi_{j_2}$. But then to avoiding an occurrence of the pattern $2134$ involving $\pi_{j_1}$ and $\pi_{j_2}$, $\pi^{\prime}$ must be the increasing permutation $n(n-1)\cdots (i_{\pi}+1)$.
\end{itemize}
 The backward direction is not difficult to see using considerations above.
\qed

Initial values for $I^{2134}_n$ are $1, 1, 3, 13, 67, 369, 2117, 12578,\ldots$ for $n\geq 1$, and this sequence is not in the OEIS~\cite{oeis}.

\begin{thm} We have
\begin{equation}\label{2134-case}
A^{2134}(x,q) = 1+ I^{2134}(x,q) + x \left(A^{2134}(x,q)-1\right) + \frac{x}{1-xq} \left(I^{213}(x,q)-x\right).\end{equation}
Moreover, for $E(x)$ and $E_n$ defined in Lemma~\ref{ge-1234}, we have
$$I^{2134}(x)=(1-x) E(x) -xC(x) -1 +3 x  + \frac{x^2}{1-x},$$
$I^{2134}_1=1$ and for $n\geq 2$,
 $$
I^{2134}_{n} = E_{n} - E_{n-1} -C_{n-1}+1.
 $$
\end{thm}

\proof
The identity \eqref{2134-case} follows from Lemma \ref{2134-decom}. Further, setting $q=1$ in \eqref{2134-case}, and applying Theorem~\ref{inde-length-3} and Lemma~\ref{ge-1234}, one has
  \begin{align*}
 I^{2134}(x) & =  (1-x) A^{2134}(x)-(1-x) - \frac{x}{1-x} I^{213}(x)   + \frac{x^2}{1-x}\\
 &=(1-x) A^{2134}(x)-(1-x) -x(C(x)-1)   + \frac{x^2}{1-x}\\
 &=(1-x) E(x) -xC(x)  +2x  + \frac{x^2}{1-x}-1.
 \end{align*}
 From this, we have the desired formula for $I^{2134}_n$. \qed

\ \\

\noindent
{\bf Pattern 1324.} We first give a description of $1324$-avoiding decomposable permutations.

\begin{lem}\label{1324-decom}
Let $\pi=\pi^{(1)}\pi^{\prime}\in \mathfrak{S}_n$, where $\pi^{(1)}$ is a permutation of $\{1,\ldots,i_{\pi}\}$. Then $\pi$ is a $1324$-avoiding decomposable permutation if and only if $\pi^{(1)}$ is $132$-avoiding and $\pi^{\prime}$ is $213$-avoiding.
\end{lem}
\proof
For the forward direction, since $\pi=\pi_1\cdots\pi_n$ is decomposable, we have $1\leq i_{\pi}\leq n-1$. Moreover,  $\pi^{(1)}$ is 132-avoiding, or else, along with $\pi_n$ an occurrence of the pattern $1324$ would be formed. Also, $\pi^{\prime}$ is 213-avoiding or else, along with $\pi_1$ an occurrence of the pattern $1324$ would be formed. The backward direction is not difficult to see, which completes the proof.
 \qed

Initial values for $I^{1324}_n$ are $1, 1, 3, 13, 69, 396, 2355, 14363, \ldots$, and this sequence is not in the OEIS~\cite{oeis}.

\begin{thm}
We have
\begin{equation}\label{1324-case}
 A^{1324}(x,q) = 1+ I^{1324}(x,q) + I^{132}(x,q) \left(A^{213}(x,q)-1\right).
\end{equation}
Also,
\begin{align*}
I^{1324}(x) = A^{1324}(x) -1 -(1-x) \left( C(x)-1\right)^2.
\end{align*}
Moreover,
$I^{1324}_1=1$ and for $n\geq 2$,
$$I^{1324}_n=A^{1324}_{n} -C_{n+1}+3C_n-2C_{n-1}.$$
\end{thm}

\proof

The identity \eqref{1324-case} follows from Lemma \ref{1324-decom}. Further,
setting $q=1$ in \eqref{1324-case} and applying Lemma \ref{catalan} and Theorem \ref{thm-indecom-132}, we obtain
 \begin{align*}
 I^{1324}(x) & =  A^{1324}(x) -1-I^{132}(x) \left(A^{213}(x)-1\right)\\
 &=A^{1324}(x) -1 -(1-x) \left( C(x)-1\right)^2\\
 &=A^{1324}(x) -1 -(1-x)\left[ \frac{C(x)-1}{x}-2C(x)+1\right].\\
\end{align*}
Hence,  $I^{1324}_1=1$ and for $n\geq 2$,
 \begin{align*}
I^{1324}_{n}& = A^{1324}_{n} - (C_{n+1}-C_n)+2(C_n-C_{n-1})\\
&= A^{1324}_{n} -C_{n+1}+3C_n-2C_{n-1}.
\end{align*}
This completes the proof.
\qed

\ \\

\noindent
{\bf Pattern 1234.} Decomposable 1234-avoiding permutations can be described as follows.

\begin{lem}\label{1234-decom}
Let $\pi = \pi^{(1)}\pi^{\prime}\in \mathfrak{S}_n$, where $\pi^{(1)}$ is a permutation of $\{1,\ldots,i_{\pi}\}$. Then $\pi$ is a $1234$-avoiding decomposable permutation  if and only if one of the following two conditions holds:
 \begin{itemize}
  \item $\pi^{(1)}=i_{\pi}(i_{\pi}-1)\cdots1$, $\pi^{\prime}$ is $123$-avoiding and $1\le i_{\pi}\le n-1$, or
  \item  $\pi^{(1)} \neq i_{\pi}(i_{\pi}-1)\cdots1$ is $123$-avoiding, $\pi^{\prime}= n(n-1)\cdots (i_{\pi}+1)$, and $3 \le i_{\pi}\le n-1$.
\end{itemize}
\end{lem}

\proof
Since $\pi$ is 1234-avoiding, then $\pi^{(1)}$ must be 123-avoiding, or else there would be an occurrence of the pattern 1234 involving an element in $\pi^{\prime}$. Thus, the longest increasing sequence in $\pi^{(1)}$ is at most of length 2. There are two cases to consider.
\begin{itemize}
\item $\pi^{(1)}=i_{\pi}(i_{\pi}-1)\cdots 1$. Then, clearly, $\pi^{\prime}$ must be $123$-avoiding.
\item The longest increasing subsequence in $\pi^{(1)}$ is exactly of length 2. But then, since
$\pi^{(1)}$ is indecomposable, we have $i_{\pi}>2$ and $\pi^{\prime}$ must be 12-avoiding, that is, $\pi^{\prime}= n(n-1)\cdots (i_{\pi}+1)$.
\end{itemize}
 This completes the proof.
 \qed

Initial values for $I_n^{1234}$ are $1, 1, 3, 13, 69, 400, 2390, 14545,\ldots$ for $n\geq 1$, and this sequence is not in the OEIS~\cite{oeis}.

\begin{thm}\label{thm-1234}
We have

$A^{1234}(x,q) = $
\begin{equation}\label{1234-case}
 1+ I^{1234}(x,q) + \frac{x}{1-xq} \left(A^{123}(x,q)-1\right)+ \frac{x}{1-xq} \left(I^{123}(x,q)-\frac{x}{1-xq}\right).
\end{equation}
Also, for $E(x)$ defined in Lemma~\ref{ge-1234},
\begin{align*}
 I^{1234}(x)
 & =  E(x) - \frac{2x}{1-x} C(x)+\frac{x^3}{(1-x)^3}+\frac{x^2}{(1-x)^2}+\frac{2 x}{1-x}-1.
 \end{align*}
 Moreover, $I^{1234}_1=1$ and for $n\geq 2$ and for $E_n$ defined in Lemma~\ref{ge-1234},
\begin{align*}
 I^{1234}_{n} & = E_{n}  -2 \sum_{i=0}^{n-1} C_i + \frac{1}{2}\left( n^2-n-4\right).
\end{align*}
\end{thm}

\proof
The identity (\ref{1234-case}) follows from Lemma \ref{1234-decom}. Further,
setting $q=1$ in \eqref{1234-case}, and applying Lemmas \ref{catalan} and \ref{ge-1234} and Theorem \ref{thm-123}, we have
 \begin{align*}
 I^{1234}(x) & =  A^{1234}(x) -1 - \frac{x}{1-x} \left(A^{123}(x)-1\right) - \frac{x}{1-x} \left(I^{123}(x)-\frac{x}{1-x}\right)\\
 & =  E(x) -1- \frac{x}{1-x} \left(2C(x)-2-\frac{x^2}{(1-x)^2}-\frac{x}{1-x}\right)\\
 & =  E(x) - \frac{2x}{1-x} C(x)+\frac{x^3}{(1-x)^3}+\frac{x^2}{(1-x)^2}+\frac{2 x}{1-x}-1.
 \end{align*}
Hence, it follows that $I^{1234}_1=1$ and for $n\geq 2$
\begin{align*}
 I^{1234}_{n} & = E_{n}  -2 \sum_{i=0}^{n-1} C_i + \frac{1}{2}\left( n^2-n+4\right)\\
 &=E_{n}  -2 \sum_{i=1}^{n-1} C_i + \frac{n(n-1)}{2}.
\end{align*}

 This completes the proof.
 \qed

\subsubsection{Pattern \texorpdfstring{$12\cdots k$}{Lg} with \texorpdfstring{$k\ge 3$}{Lg}}

Here  we consider patterns of the form $12\cdots k$, where $k\geq 3$, which generalizes our considerations for patterns $123$ and $1234$. First, we give a description of $12\cdots k$-avoiding decomposable permutations in the following lemma, whose proof is trivial and thus is omitted.

\begin{lem}\label{aaa}
If $\pi=\pi^{(1)}\pi^{\prime}$ is a $12\cdots k$-avoiding decomposable permutation of length $n$, where $\pi^{(1)}$ is a permutation of $\{1,\ldots,i_{\pi}\}$.
Then there exists $m$, $1 \le m\le k-2$, such that
the longest increasing subsequence in $\pi^{(1)}$ is exactly of length $m$ and  $\pi^{\prime}$ is $12\cdots (k-m)$-avoiding.
\end{lem}

Now we can enumerate $12\cdots k$-avoiding indecomposable permutations.
\begin{cor}
We have
\begin{align*}
I^{12\cdots k}(x,q)=A^{12\cdots k}(x,q)-\sum_{m=1}^{k-2}\left(I^{12\cdots (m+1)}(x,q)-I^{12\cdots m}(x,q)\right) \cdot \left(A^{12\cdots(k-m)}(x,q)-1\right)-1.
\end{align*}
\end{cor}

\proof
By Lemma \ref{aaa}, we have
$$
A^{12\cdots k}(x,q)=1+I^{12\cdots k}(x,q)+\sum_{m=1}^{k-2}\left(I^{12\cdots (m+1)}(x,q)-I^{12\cdots m}(x,q)\right) \cdot \left(A^{12\cdots(k-m)}(x,q)-1\right),
$$
from which the result follows.
\qed

For example, when $k=3$, we have
\begin{align*}
A^{123}(x) & = 1+I^{123}(x)+\left(I^{12}(x)-I^{1}(x)\right) \cdot \left(A^{12}(x)-1\right)\\
& = 1+I^{123}(x)+\frac{x}{1-x} \cdot \frac{x}{1-x},
\end{align*}
and hence
\begin{align*}
I^{123}(x) = A^{123}(x) - 1-\frac{x^2}{(1-x)^2} = \frac{1-\sqrt{1-4 x}}{2 x} - 1-\frac{x^2}{(1-x)^2},
\end{align*}
which coincides with Theorem \ref{thm-123}.
Note that we used the facts that  $I^{(1)}(x)=0$ and $A^{12}(x)-1=I^{(12)}(x)=x+x^2+\cdots=\frac{x}{1-x}$.

When $k=4$, we have
$$
A^{1234}(x) = 1 + I^{1234}(x) + \left(I^{12}(x)-I^{1}(x)\right) \cdot \left(A^{123}(x)-1\right) +\left(I^{123}(x)-I^{12}(x)\right) \cdot \left(A^{12}(x)-1\right),
$$
and hence
\begin{align*}
 I^{1234}(x)
 & =  E(x) - \frac{x}{1-x} \left(2C(x)-2-\frac{x^2}{(1-x)^2}-\frac{x}{1-x}\right)-1,
 \end{align*}
 which coincides with Theorem \ref{thm-1234}. Note that we used the fact that $A^{123}(x)=C(x).$

\section{Indecomposable permutations avoiding vincular non-consecutive patterns of length 3}\label{vin-sec}

For a pattern of the form $a\underline{bc}$, its reverse complement gives a pattern of the form $\underline{xy}z$. Thus, since the composition of reverse and complement preserves the property of being indecomposable,  we only need to consider six cases of  vincular patterns of length 3, which are
$1\underline{23}$, $1\underline{32}$, $2\underline{13}$, $2\underline{31}$, $3\underline{12}$, and $3\underline{21}$. Two of these cases can be reduced to classical pattern-avoidance.

Indeed, it was shown in \cite{Claesson2001Generalized} that a permutation avoids the pattern $2\underline{13}$ if and only if it avoids the pattern 213. Applying the complement operation, this implies that a permutation avoids $2\underline{31}$ if and only if it avoids 231. Thus, $I^{2\underline{13}}_n=I^{213}_n$ and $I^{2\underline{31}}_n=I^{231}_n$ and Theorems~\ref{thm-indecom-132} and \ref{inde-length-3} can be applied, respectively.\\

\noindent
{\bf Pattern $1\underline{23}$.} We first give a description of $1\underline{23}$-avoiding decomposable permutations.

\begin{lem}\label{lem-123}
Let $\pi=\pi^{(1)}\pi^{\prime}\in \mathfrak{S}_n$, where $\pi^{(1)}$ is a permutation of $\{1,\ldots,i_{\pi}\}$. Then $\pi=\pi_1\cdots \pi_n$ is a $1\underline{23}$-avoiding decomposable permutation if and only if $\pi^{(1)}=\pi_1 \cdots \pi_{i_{\pi}-1}1$
and $\pi^{\prime}=n(n-1)\cdots (i_{\pi}+1)$ for $1\leq i_{\pi}\leq n-1$, where
$\pi_1 \cdots \pi_{i_{\pi}-1}$ is a $1\underline{23}$-avoiding permutation of $\{2,3,\ldots,i_{\pi}\}$.
\end{lem}

\proof
Since $\pi$ is decomposable, we have $1\leq i_{\pi}\leq n-1$. It is clear that $\pi^{(1)}$ is a $1\underline{23}$-avoiding indecomposable permutation. We claim that $\pi_{i_{\pi}}=1$, since otherwise 1,  $\pi_{i_{\pi}}$ and $\pi_{i_{\pi}+1}$ will form the pattern $1\underline{23}$. Further, clearly $\pi_{i_{\pi}+1}> \pi_{i_{\pi}+2}>\cdots>\pi_n$, or else there would be an occurrence of the pattern $1\underline{23}$ involving 1.

On the other hand, it is easy to see that if $\pi^{(1)}$ and $\pi^{\prime}$ satisfy the conditions then $\pi$ is $1\underline{23}$-avoiding.  This completes the proof.
\qed

Initial values for $I_n^{1\underline{23}}$ for $n\geq 1$ are $1, 1, 3, 11, 43, 179, 801, \ldots$ and this sequence is not in the OEIS~\cite{oeis}.

\begin{thm}
We have

\begin{equation}\label{1-23-case}
 A^{1\underline{23}}(x,q)=1+I^{1\underline{23}}(x,q)+\frac{x}{1-xq}
\left[\left(A^{1\underline{23}}(x,q)-1\right)xq+x\right].
\end{equation}
Also,
$$I^{1\underline{23}}(x)=\frac{1-x-x^2}{1-x}B(x)-1,$$
where $B(x)$ is the generating function for the Bell numbers.

Moreover, $I^{1\underline{23}}_1=1$ and for $n\geq 2$,
$$I^{1\underline{23}}_n=B_n-\sum_{i=0}^{n-2}B_i.$$
\end{thm}
\proof
By Lemma \ref{lem-123}, the generating function for $1\underline{23}$-avoiding decomposable permutations is
$\frac{x}{1-xq}\left[\left(A^{1\underline{23}}(x,q)-1\right)xq+x\right]$, hence (\ref{1-23-case}) follows.

Letting $q=1$ in (\ref{1-23-case}), we have
$$I^{1\underline{23}}(x)=\left(1-\frac{x^2}{1-x}\right)A^{1\underline{23}}(x)-1.$$
Combing with Lemma \ref{bell}, we obtain that
$$I^{1\underline{23}}(x)=\frac{1-x-x^2}{1-x}B(x)-1.$$
Together with the fact that
$\frac{1}{1-x}=\sum_{i\geq 0}x^i$, we have $I^{1\underline{23}}_1=1$ and for $n\geq 2$,
$$I^{1\underline{23}}_n=B_n-\sum_{i=0}^{n-2}B_i.$$
 This completes the proof.
 \qed

\ \\

\noindent
{\bf Pattern $1\underline{32}$.} We first give a description of $1\underline{32}$-avoiding decomposable permutations.

\begin{lem}\label{lem-132}
Let $\pi=\pi^{(1)}\pi^{\prime}\in \mathfrak{S}_n$, where $\pi^{(1)}$ is a permutation of $\{1,\ldots,i_{\pi}\}$. Then $\pi=\pi_1\cdots \pi_n$ is a $1\underline{32}$-avoiding decomposable permutation if anf only if $\pi^{(1)}$
is $1\underline{32}$-avoiding and $\pi^{\prime}=i_{\pi}(i_{\pi}+1)\cdots n$ for $1\leq i_{\pi}\leq n-1$.
\end{lem}
\proof
If $\pi$ is $1\underline{32}$-avoiding, then clearly $\pi^{(1)}$ and $\pi^{\prime}$ are both $1\underline{32}$-avoiding. Moreover, we must have
$\pi_{i_{\pi}+1}<\pi_{i_{\pi}+2}<\cdots<\pi_n$, or else there would be an occurrence of the pattern $1\underline{32}$ involving 1.

On the other hand, it is clear that if $\pi^{(1)}$ and $\pi^{\prime}$ satisfy the given conditions then $\pi$ is $1\underline{32}$-avoiding, which completes the proof.
 \qed

Initial values for $I_n^{1\underline{32}}$ for $n\geq 1$ are  $1, 1, 3, 10, 37, 151, 674,\ldots$, which are essentially the sequences A005493 and A138378 in the OEIS~\cite{oeis} that have several combinatorial interpretations. In particular, this sequence counts $1\underline{32}$-avoiding permutations that end with a rise, that is, with an occurrence of the pattern $\underline{12}$, which leads us to the following theorem.

\begin{thm}\label{thm-bij} For $n\geq 2$, the number of $1\underline{32}$-avoiding  indecomposable permutations in $\mathfrak{S}_n$ is equal to that of $1\underline{32}$-avoiding permutations in $\mathfrak{S}_n$ that end with a rise. \end{thm}

\proof
Let $I$ and $R$ be the first and the second sets, respectively, in the statement of the theorem. We provide a recursive bijection $f$ from $I$ to $R$ proving the theorem with the base case $f(21)= 12$.

The set of $1\underline{32}$-avoiding permutations can be subdivided into three disjoint subsets:
\begin{itemize}
\item $S_1$, all $1\underline{32}$-avoiding permutations ending with 1;
\item $S_2$, all $1\underline{32}$-avoiding  permutations ending with $n$;
\item $S_3$, all other $1\underline{32}$-avoiding permutations.
\end{itemize}
It is straightforward to see that to the right of 1 in a $1\underline{32}$-avoiding permutation we must have an increasing order of elements. But then in $S_3$, $n$ must be to the left of 1. Thus, a permutation in $S_3$ belongs to both $I$ and $R$  and we map it to itself.  Further, it is easy to see that $S_1$ is a subset of $I$ but it is disjoint from $R$, while $S_2$ is a subset of $R$ but it is disjoint from $I$. For a permutation $\pi=\pi_1\cdots\pi_{n-1}1\in S_1$ we define its image recursively as
$$f(\pi)=f((\pi_1-1)\cdots(\pi_{n-1}-1))n.$$
The map $f$ described by us is easy to see to be a bijection.
\qed

Next we enumerate $1\underline{32}$-avoiding indecomposable permutations.

\begin{thm}
We have
\begin{equation}\label{1-32-case}
A^{1\underline{32}}(x,q)=1+I^{1\underline{32}}(x,q)+\frac{x}{1-x}I^{1\underline{32}}(x,q).
\end{equation}
Also,
$$I^{1\underline{32}}(x)=(1-x)B(x)-1,$$
where $B(x)$ is the generating function for the Bell numbers.
Moreover, $I^{1\underline{32}}_n=1$ and for $n\geq 2$
$$I^{1\underline{32}}_n=B_n-B_{n-1}.$$ \end{thm}
\proof

By Lemma \ref{lem-132}, the generating function for $1\underline{32}$-avoiding decomposable permutations is $\frac{x}{1-x}I^{1\underline{32}}(x,q)$ from which (\ref{1-32-case}) follows.

Letting $q=1$ in (\ref{1-32-case}) and combing with Lemma \ref{bell}, we obtain that
$$I^{1\underline{32}}(x)=(1-x)A^{1\underline{32}}(x)-1=(1-x)B(x)-1.$$
Hence, it follows that
$$I^{1\underline{32}}_n=B_n-B_{n-1}$$ for $n\geq 2$.
 This completes the proof.\qed

\ \\

\noindent
{\bf Patterns $3\underline{12}$ and $3\underline{21}$.} We first give a description of $3\underline{12}$-avoiding and $3\underline{21}$-avoiding decomposable permutations.

 \begin{lem}
 Let $\pi=\pi^{(1)}\pi^{\prime}\in \mathfrak{S}_n$, where $\pi^{(1)}$ is a permutation of $\{1,\ldots,i_{\pi}\}$. Then  $\pi$ is $3\underline{12}$- (resp., $3\underline{21}$-)avoiding if and only if $\pi^{(1)}$ and $\pi^{(2)}$ are both $3\underline{12}$- (resp., $3\underline{21}$-)avoiding.
 \end{lem}

Initial values for $I_n^{3\underline{12}}=I_n^{3\underline{21}}$ for $n\geq 1$ are $1, 1, 2, 6, 22, 92, 426,\ldots$, and this is the sequence A074664 in the OEIS~\cite{oeis} that have several combinatorial interpretations.  In particular, this sequence counts the number of {\em irreducible set  partitions} of $[n]$,
which can be easily seen from the bijections in \cite{Claesson2001Generalized}.
For more information on irreducible set  partitions, see \cite{Chen2011bijection}.

 \begin{thm}
 We have
 $$I^{3\underline{12}}(x)=I^{3\underline{21}}(x)=1-\frac{1}{B(x)}.$$
 \end{thm}

 \proof
Since 3\underline{12} is an irreducible pattern, Lemma~\ref{inde} can be applied to obtain
\begin{align*}
 A^{3\underline{12}}(x) = \frac{1}{1-I^{3\underline{12}}(x)}.
\end{align*}
The desired result now follows from Lemma~\ref{bell}.\qed

\section{Concluding remarks}\label{research}

The notion of indecomposable permutations proved to be useful in various contexts, e.g. in obtaining non-trivial enumeration and equidistribution results on permutations \cite{CKS}.

In this paper, we gave a compete classification of indecomposable permutations avoiding a classical pattern of length 3 or 4, and of indecomposable permutations avoiding a non-consecutive vincular pattern of length 3. Also, we provided a recursive formula for enumerating $12\cdots k$-avoiding indecomposable permutations for $k\geq 3$. The descent statistic is taken into account in several of our results.

A natural direction of further research is in extending our studies of indecomposable permutations to other patterns, e.g. vincular patterns of length 4. Also, one can look at avoiding more than one pattern at the same time. Other statistics can be included in enumerative results.

Finally, one can establish a number of bijective results linking pattern avoiding indecomposable permutations to other structures (Theorem~\ref{thm-bij} is one such example). For instance, the sequence A005493 in the OEIS \cite{oeis} has many interesting combinatorial interpretations that one could try to link in a bijective way to $1\underline{32}$-avoiding indecomposable permutations.

\section*{\bf Acknowledgments}
The work of the first and the third authors was supported by the 973 Project, the PCSIRT Project of the Ministry of Education and the National Science Foundation of China. The second author is grateful to the administration of the Center for Combinatorics at Nankai University for their hospitality during the author's stay in November -- December 2015. Also, the authors are grateful to Filippo Disanto for brining to out attention the paper \cite{Disanto}.

%\bibliographystyle{abbrv}
%\bibliography{../../../ref}

% \begin{thebibliography}{10}
%
% \bibitem{CK} A. Claesson, S. Kitaev. Classification of
% bijections between 321- and 132-avoiding permutations, {\em
% S\'eminaire Lotharingien de Combinatoire} {\bf B60d} (2008), 30 pp.
%
%
% \bibitem{Kitaev2011Patterns}
% S.~Kitaev. {\em Patterns in permutations and words}, Springer, 2011.
%
% \end{thebibliography}

\end{document}